\documentclass[11pt]{article}

\usepackage{amssymb}
\usepackage{amsmath}
\usepackage{theorem}
\usepackage{epsfig}
\usepackage{verbatim}
\usepackage{graphicx}

\usepackage{color}

\newtheorem{thm}{Theorem}[section]

\newtheorem{lemma}[thm]{Lemma}

\newcommand{\R}{\Bbb{R}}
\newcommand{\C}{\Bbb{C}}

\newcommand{\N}{\Bbb{N}}
\newcommand{\T}{\mathbb{T}}
\newcommand{\D}{\displaystyle}

\newcommand{\grad}{\nabla}

\newcommand{\gradp}{\grad^{\bot}}

\newcommand{\da}{\partial_\alpha}

\newcommand{\la}{\Lambda}

\newcommand{\sign}{{\rm sign}\thinspace}

\newcommand{\al}{\alpha}
\newcommand{\ep}{\varepsilon}

\newcommand{\pa}{\partial}

\newcommand{\be}{\beta}

\begin{document}

\title{Turning waves and breakdown for incompressible flows}
\author{\'Angel Castro, Diego C\'ordoba, Charles Fefferman,\\
 Francisco Gancedo and Mar\'ia L\'opez-Fern\'andez.}
\date{November 26, 2010}

\maketitle

\begin{abstract}

We consider the evolution of an interface generated between two
immiscible incompressible and irrotational fluids. Specifically we
study the Muskat and water wave problems. We show that starting
with a family of initial data given by $(\al,f_0(\al))$, the
interface reaches a regime in finite time in which is no longer a
graph. Therefore there exists a time $t^*$ where the solution of the
free boundary problem parameterized as $(\al,f(\al,t))$ blows-up:
$\|\da  f\|_{L^\infty}(t^*)=\infty$. In particular, for the Muskat
problem, this result allows us to reach an unstable regime, for
which the Rayleigh-Taylor condition changes sign and the solution breaks down.

\end{abstract}

\maketitle


\section{Introduction}

Here we study two problems of Fluids Mechanics concerning the
evolution of two incompressible fluids of different characteristics
in 2D. We consider that both fluids are immiscible and of different
constant densities $\rho^1$ and $\rho^2$, modeling the dynamics of
an interface that separates the domains $\Omega^1(t)$ and
$\Omega^2(t)$. That is, the liquid density $\rho=\rho(x,t)$,
$(x,t)\in \R^2\times\R^+$, is defined by
\begin{equation}\label{density}
\rho(x,t)=\left\{\begin{array}{cl}
                    \rho^1,& x\in\Omega^1(t)\\
                    \rho^2,& x\in\Omega^2(t)=\mathbb{R}^2 - \Omega^1(t),
                 \end{array}\right.
\end{equation}
and satisfies the conservation of mass equation
\begin{align}
\begin{split}\label{mass}
\rho_t+v\cdot\grad\rho&=0, \\
\grad\cdot v&=0,
\end{split}
\end{align}
where $v = (v_1(x,t), v_2(x,t))$ is the velocity field.
 With a free boundary parameterized by
$$
\partial \Omega^j(t)=\{z(\al,t)=(z_1(\al,t),z_2(\al,t)):\al\in\R\},
$$
we consider open curves vanishing at infinity
$$\D\lim_{\al\rightarrow\infty}(z(\al,t)-(\al,0))=0,$$
or periodic in the space variable
\begin{equation*}
z(\al+2k\pi,t)=z(\al,t)+2k\pi(1,0).
\end{equation*}
The scalar vorticity, $\gradp\cdot v$, has the form
\begin{equation}\label{vorticity}
\gradp\cdot v(x,t)=\omega(\al,t)\delta(x-z(\al,t)),
\end{equation}
i.e. the vorticity is a Dirac measure on $z$ defined by
$$
<\gradp\cdot v,\eta>=\int_{\R}\omega(\al,t)\eta(z(\al,t))d\al,
$$
with $\eta(x)$ a test function. The system is closed by using one of
the following fundamental fluid motion equations:\\

Darcy's law
\begin{equation}\label{Darcy}
\frac{\mu}{\kappa}v=-\grad p-g\rho(0,1),
\end{equation}
or\\

Euler equations
\begin{equation}\label{Euler}
\rho(v_t+v\cdot\grad v)=-\grad p-g\rho(0,1).
\end{equation}
Here $p=p(x,t)$ is pressure, $g$ gravity, $\mu$ viscosity and
$\kappa$ permeability of the isotropic medium.

The Muskat problem \cite{Muskat} is given by
(\ref{density},\ref{mass},\ref{Darcy}) which considers the dynamics
of two incompressible fluids of different densities throughout
porous media and Hele-Shaw cells \cite{S-T,Howison}. In this last
setting the fluid is trapped between two fixed parallel plates, that
are close enough together, so that the fluid essentially only moves
in two directions \cite{H-S}.

Taking $\rho^1=0$, equations
(\ref{density},\ref{mass},\ref{vorticity},\ref{Euler}) are known as
the water waves problem (see \cite{BL} and references therein),
modeling the dynamics of the contour between an inviscid fluid with
density $\rho^2$ and vacuum (or air) under the influence of gravity.

Condition \eqref{vorticity} (deduced by \eqref{Darcy}, assumed for
\eqref{Euler}) allows us to write the evolution equation in terms of
the free boundary as follows. One could recover the velocity field
from \eqref{vorticity} by means of Biot-Savart law
\begin{equation*}\label{BS}
v(x,t)=\gradp\Delta^{-1}(\gradp\cdot v)(x,t)=\frac{1}{2\pi}\int_\R
\frac{(x-z(\al,t))^{\bot}}{|x-z(\al,t)|^2}\omega(\al,t)d\al,
\end{equation*}
applying the delta measure with amplitude $\omega$. Taking limits on
the above equation approaching the boundary in the normal direction
inside $\Omega^j$, the velocity is shown to be discontinuous in the
tangential direction, but continuous in the normal, and given by the
Birkhoff-Rott integral of the amplitude $\omega$ along the interface
curve:
\begin{equation*}\label{BR}
BR(z,\omega)(\al,t)=\frac{1}{2\pi}PV\int_\R
\frac{(z(\al,t)-z(\beta,t))^\bot}{|z(\al,t)-z(\beta,t)|^2}\omega(\beta,t)
d\beta.
\end{equation*}
Above $PV$ denotes principal value. It yields the curve velocity
from which one can subtract any term $c$ in the tangential without
modifying the geometry of the interface
\begin{equation}\label{em}
z_t(\al,t)=BR(z,\omega)(\al,t)+c(\al,t)\da z(\al,t).
\end{equation}
Understanding the problem as weak solutions of
(\ref{density},\ref{mass},\ref{Darcy})  or
(\ref{density},\ref{mass},\ref{vorticity},\ref{Euler}), the
continuity of the pressure on the free boundary follows. Therefore,
taking limits in Darcy's law from both sides and subtracting the
results in the tangential direction, it is easy to close the system
for Muskat (in this paper we consider two fluids with the same viscosity):
\begin{equation}\label{cDarcy}
\omega(\al,t)=-(\rho^2-\rho^1)\frac{\kappa g}{\mu}\da z_2(\al,t).
\end{equation}
In a similar way for water waves, Euler equations yield
\begin{align}
\begin{split}\label{cEuler}
\omega_t(\al,t)&=-2\partial_t BR(z,\omega)(\al,t)\cdot
\da z(\al,t)-\da( \frac{|\omega|^2}{4|\da z|^2})(\al,t) +\da (c\, \omega)(\al,t)\\
&\quad +2c(\al,t)\da BR(z,\omega)(\al,t)\cdot\da z(\al,t)+2g\da
z_2(\al,t).
\end{split}
\end{align}
Then, the two contour equations are set by (\ref{em},\ref{cDarcy})
and (\ref{em},\ref{cEuler}).

 For these models the well-posedness turns out to be false for some settings.
Rayleigh \cite{Ray} and Saffman-Taylor \cite{S-T} gave a condition that must
be satisfied for the linearized model in order to exist a solution
locally in time: the normal component of the pressure
gradient jump at the interface has to have a distinguished sign.
This is known as the Rayleigh-Taylor condition. It reads
$$\sigma(\al,t) = -(\nabla p^2(z(\al,t),t) - \nabla
p^1(z(\al,t),t))\cdot\partial^{\perp}_{\al}z(\al,t)>0,$$ where
$\nabla p^j(z(\al,t),t)$ denotes the limit gradient of the pressure
obtained approaching the boundary in the normal direction inside
$\Omega^j(t)$.

  An easy linearization around a flat contour
$(\al,f(\al,t))$, allows us to find
$$
f_t=\frac{1}{2}H(\omega)
$$
where $H$ is the Hilbert transform which symbol on the Fourier side
is given by $\widehat{H}=-i\, \sign (\xi)$. The equations
\begin{align*}
\omega&=-(\rho^2-\rho^1)\frac{\kappa g}{\mu}\da f,&
\mbox{(linear Muskat)}\\
 \omega_t&=2g\da f,&\mbox{(linear water waves)}
\end{align*} show the parabolicity of the Muskat problem when the denser fluid is below
($\rho^2>\rho^1$) and the dispersive character of water waves.

1. There is a wide literature on the Muskat problem and the dynamics
of two fluids in a Hele-Shaw cell. There are works considering the
case of a viscosity jump neglecting the effect of gravity  \cite{SCH,Yi}.
Local-existence in a more general situation (with discontinuous
viscosity and density) is shown in \cite{ADP} and also treated in
\cite{Am}.  A different approach to prove local-existence can be
found in \cite{DY} for the setting we are considering in this paper.
The Rayleigh-Taylor stability depends upon the sign of
$(\rho^2-\rho^1)\da z_1(\al,t)$ indicating that the heavier fluid
has to be below in the stable case. If the lighter fluid is below,
the problem have been shown to be ill-posed \cite{DY}.
Global-existence results for small initial data can be found in
\cite{Peter,Yi2,SCH,DY,Esch2}. For large initial curves and
parameterized by $(\al,f(\al,t)),$ there are maximum principles for
the $L^\infty$ and $L^2$ norms of $f$, and decay rates, together with
global-existence for Lipschitz curves if $\|\partial_\al
f\|_{L^\infty}(0)<1$ \cite{DP2, ccgs}.

2. The water waves problem have been extensively considered (see
\cite{ADP2,BL} and references therein). For sufficiently smooth free
boundary, the Rayleigh-Taylor condition remains positive with no
bottom considerations \cite{Wu}, a fact that was used to prove
local-existence \cite{Wu}. The Rayleigh-Taylor stability can play a
different role for the case of non "almost" flat bottom
\cite{Lannes}. Recently, for small initial data, exponential time of
existence has been proven in two dimensions \cite{Wu3}, and
global-existence in the three dimensional case (two dimensional
interface) \cite{GMS,Wu4}.

The results that we announce below will be given in the forthcoming papers \cite{ADCPM} and \cite{ADCPM2}.

\section{Rayleigh-Taylor breakdown for Muskat}

 This section is devoted to show the main ingredients to prove
the Theorem \ref{perturbativo} below. We consider the function
\begin{equation}\label{arccord}
F(z_0)(\al,\beta)=\frac{|\beta|^2}{|z_0(\al)-z_0(\al-\beta)|^2},
\end{equation}
if $F(z_0)\in L^\infty(\R^2)$ then the curve $z_0$ satisfies the arc-chord condition. 
We say that the Rayleigh-Taylor of the solution of the Muskat problem breaks down in finite time if for initial data $z_0$ satisfying $\sigma(\al,0)=(\rho^2-\rho^1)\da z_1(\al,0)>0$ there exists a time $t^*>0$ for which $\sigma(\al,t^*)$ is strictly
negative in a nonempty open interval.

\begin{thm}\label{perturbativo}
There exists a non-empty open set of initial data in $H^4$, satisfying
Rayleigh-Taylor and arc-chord conditions, for which the
Rayleigh-Taylor condition of the solution of the Muskat problem
(\ref{density}, \ref{mass}, \ref{Darcy}) breaks down in finite time.
\end{thm}

 After choosing the appropriate tangential term and a integration by parts, the contour equation reads
\begin{equation}\label{ec1d}
z_t(\al,t) =\frac{\rho^2-\rho^1}{2\pi}PV\int
\frac{(z_1(\al,t)-z_1(\beta,t))}{|z(\al,t)-z(\beta,t)|^2}(\da
z(\al,t) - \da z(\beta,t)) d\beta.
\end{equation}

The steps of the proof:\\

1. First, for any initial curve $z_0(\al)=z(\al,0)$ in $H^4$ that
satisfy R-T
$$(\rho^2-\rho^1)\da z_1(\al,0)>0$$
and the arc-chord condition
then the solution to the Muskat problem $z(\al,t)$ becomes analytic
for $0<t<T$. Moreover, $z(\al,t)$ is real analytic in a strip
$$S(t)=\{\al+i\zeta:|\zeta|<ct\}$$ for $t\in (0,T)$ where $c$
depends only on
$$\inf(0)=\inf_{\al}\frac{\da z_1(\al,0)}{|\da z(\al,0)|^2}.$$
The proof follows by controlling the quantities extended on $S(t)$:
$$
f(\al+i\zeta,t)=\frac{\da z_1(\al+i\zeta,t)}{|\da
z(\al+i\zeta,t)|^2},
$$
$F(z)(\al+i\zeta,\beta,t)$ by using \eqref{arccord} and norms
$$
\|F(z)\|_{L^\infty(S)}(t)=\sup_{\al+i\zeta\in
S(t),\beta\in\T}|F(z)(\al+i\zeta,\beta)|,
$$
$$
\|z\|^2_{L^2(S)}(t)=\sum_{\pm}\int_{\T}|z(\al\pm ict,t)|^2 d\al,
$$
$$
\|z\|^2_{H^j(S)}(t)=\|z\|^2_{L^2(S)}(t)+\sum_{\pm}\int_{\T}|\da^jz(\al\pm
ict,t)|^2 d\al,
$$
for $j\in\N$,
$$
\inf(t)=\inf_{\al+i\zeta\in S(t)}\Re(f)(\al+i\zeta).
$$
Then the quantity
$$\|z\|^2_{RT}(t)=\|z\|^2_{H^4(S)}(t)+\|F(z)\|_{L^\infty(S)}(t)+1/(\inf(t)-c-K\|\Im(f)\|_{H^2(S)}(t))$$
satisfies
$$
\frac{d}{dt}\|z\|_{RT}(t)\leq C\|z\|^k_{RT}(t),
$$
for $C$, $K$ and $k$ universal constants. It yields
$$
\|z\|_{RT}(t)\leq \frac{\|z\|_{RT}(0)}{(1-C\|z\|_{RT}^k(0)t)^{1/k}},
$$
providing control of the analyticity and $T=1/(C\|z\|_{RT}^k(0))$.

2. Second, there is a lower bound on the strip of analyticity, which
does not collapse to the real axis as long as the Rayleigh-Taylor is
greater than or equal to 0. Then there is a time $T$ and a solution
of the Muskat problem $z(\al,t)$ defined for $0<t\leq T$ that
continues analytically  into a complex strip if $(\rho^2-\rho^1)\da
z_1\geq 0$, where $T$ is either a small constant or it is the first
time a vertical tangent appears, whichever occurs first. We
redefine the strip
$$S(t)=\{\al+i\zeta:|\zeta|<h(t), 0<h(0)\},$$
and the quantity $\|z\|^2_{S}=\|z\|^2_{H^4(S)}+\|F(z)\|_{L^\infty(S)}$
with this new $S(t)$. For a $h(t)$ decreasing (the expression of
$h(t)$ is chosen later), we consider the evolution of the most
singular quantity
$$\sum_{\pm}\int |\da^4 z(\al\pm ih(t),t)|^2d\alpha.$$ Taking a
derivative in $t$ one finds
$$\frac{d}{dt}\sum_{\pm}\int |\da^4 z(\al\pm ih(t))|^2 d\alpha\leq \frac{h'(t)}{10}\sum_{\pm}\int \la(\da^4z)(\al\pm ih(t))\cdot\overline{\da^4 z}(\al\pm ih(t))d\al$$
$$-10 h'(t)\int\la(\da^4 z)(\al)\cdot\overline{\da^4z}(\al)d\al+2\sum_{\pm}\Re\int
\da^4 z_t(\al\pm i h(t))\cdot\overline{\da^4z}(\al \pm
ih(t))d\alpha.$$ Estimating in a wise way one obtains
$$
\frac{d}{dt}\sum_{\pm}\int|\da^4 z(\al\pm ih(t))|^2 d\alpha\leq
C\|z\|^k_{S}(t)-10 h'(t)\int\la(\da^4
z)(\al)\cdot\overline{\da^4z}(\al)d\al$$
$$
+(C\|z\|^k_{S}(t) h(t)+ \frac{1}{10}h'(t))\int \la(\da^4z)(\al\pm
ih(t))\cdot\overline{\da^4 z}(\al\pm ih(t))d\al.
$$
Therefore, choosing
$$h(t)=h(0)\exp(-10C\int_0^t\|z\|^k_{S}(r)dr)$$
eliminates the most dangerous term. The other terms are easily
controlled, giving finally
$$
\frac{d}{dt}\sum_{\pm}\int|\da^4 z(\al\pm ih(t))|^2 d\alpha\leq
C\|z\|^k_{S}(t).
$$
This allows us to reach a regime for which the boundary $z$ develops
a vertical tangent at time $T$.

3. Third, it is shown the existence of a large class of analytic
curves for which there exist a point where the tangent vector is
vertical and the velocity indicates that the curve is going to turn
up and reach the unstable regime. That is
$$
\begin{array}{ll}
a.\,\,\da z_1(\alpha)>0\mbox{ if }\alpha \neq 0,\qquad & b.\,\,\da z_1(0)=0,\\
&\\
c.\,\,\da z_2(0)>0, & d.\,\, \da v_1(0)<0,
\end{array}
$$
for analytic functions $z_1(\al)$ and $z_2(\al)$ such that $z(\al)$
satisfies the arc-chord condition. Here we consider an open curve
vanishing at infinity (being analogous in the periodic case). We
assume that $z(\alpha)$ is a smooth odd curve satisfying the
properties $a$, $b$ and $c$. Differentiating the expression
\eqref{ec1d} for the horizontal component of the velocity it is easy
to obtain
$$(\da v_1)(\al)=\int_{-\infty}^\infty \frac{(\da z_1(\al)-\da z_1(\al\!-\!\be))^2+(z_1(\al)-z_1(\al\!-\!\be))(\da^2 z_1(\al)-\da^2 z_1(\al\!-\!\be))}{|z(\al)-z(\al\!-\!\be))|^2}d\be$$
$$-2\int_{-\infty}^\infty (z_1(\al)-z_1(\al\!-\!\be))(\da z_1(\al)-\da z_1(\al\!-\!\be))\frac{(z(\al)-z(\al\!-\!\be))\cdot(\da z(\al)-\da z(\al\!-\!\be))}{|z(\al)-z(\al\!-\!\be)|^4}d\be. $$
At $\al=0$ it yields
$$(\da v_1)(0)=\int_{-\infty}^\infty\frac{(\da z_1(\be))^2+z_1(\be)\da^2 z_1(\be)}{|z(\be)|^2}d\be$$
$$-2\int_{-\infty}^\infty z_1(\be)\da z_1(\be)\frac{z_1(\be)\da z_1(\be)-z_2(\be)(\da z_2(0)-\da z_2(\be))}{|z(\be)|^4}d\be.$$
Integration by parts provides
$$\int_{-\infty}^\infty\frac{z_1(\be)\da^2 z_1(\be)}{|z(\be)|^2}d\be=-\int_{-\infty}^\infty\frac{(\da z_1(\be))^2}{|z(\be)|^2}d\be$$
$$+2\int_{-\infty}^\infty z_1(\be)\da z_1(\be)\frac{z_1(\be)\da z_1(\be)+z_2(\be)\da z_2(\be)}{|z(\be)|^4}d\be.$$
Therefore it is easy to obtain that
\begin{equation}\label{reducida}
(\da v_1)(0)=4\da z_2(0)\int_0^{\infty}\frac{z_1(\be)z_2(\be)}{|
z(\be)|^4}\da z_1(\be)d\be.
\end{equation}
Expression (\ref{reducida}) allows us to determine the sign of $(\da
v_1)(0)$. One could take
$$z_1(\be)=\frac{\be^3}{(1+\be^2)},$$
and construct the function $z_2(\be)$ in the following way: Let
$\be_1$, $\be_2$ and $\be_3$ be real increasing numbers. We pick
$z_2(\be)<c<0$ for $\be_2< \be<\infty$ and $z_2^*(\be)$ a smooth
function with the following properties
$$
\begin{array}{ll}
a.\,\,z_2^*(\be)\mbox{ is odd,}\qquad & b.\,\,(\partial_\be z_2^*)(0)>0,\\
&\\
c.\,\,z_2^*(\be)>0 \mbox{ if } \be\in (0,\be_1), & d.\,\,
z_2^*(\be)<0 \mbox{ if }\be\in (\be_1,\be_2].
\end{array}
$$
For $z_2(\be)=bz_2^*(\be)$, $0\leq \be\leq \be_2$ and $b>0$, the
velocity satisfies
$$(\da v_1)(0)< 4(\da z_2)(0)\left (\int_0^{\be_1}\frac{z_1(\be)z_2(\be)}{|
z(\be)|^4} \pa_\al
z_1(\be)d\be+\int_{\be_3}^\infty\frac{z_1(\be)z_2(\be)}{| z(\be)|^4}
\pa_\al z_1(\be)d\be\right)$$
$$=4(\da z_2)(0)\left(\int_0^{\be_1}\frac{z_1(\be)bz^*_2(\be)}{(
z_1(\beta)^2+b^2z_2^*(\beta)^2)^2}\pa_\al z_1(\be) d\be+ A\right),$$
where $A<0$. The constant $b$ large enough yields $(\da v_1)(0)<0.$
Rectifying the curve on the interval $[\be_2,\be_3]$ it is easy to obtain a smooth curve. Finally, convolving with the heat kernel the vertical component, the curve $z(\alpha)$ is approximated by an analytic one.

4. Fourth, with the initial data found in 3. and no assumption on the R-T condition, we use a modification of Cauchy-Kowalewski theorems \cite{Nirenberg,Nishida} to show that there exists an analytic solution for the Muskat problem in some interval $[-T,T]$ for a small enough $T>0$. Here we are forced to change
substantially the method in \cite{SSBF} because in this case the curve can not be parameterized as a graph, so we have to deal with the arc-chord condition. Then, with $\{X_r\}_{r>0}$ a scale of Banach spaces given by real functions that can be extended analytically on the complex
strip $S_r=\{\al+i\zeta\in\C: |\zeta|< r\}$ with norm
$$
\|f\|_r=\sum_{\pm}\int|f(\al\pm ir)|^2d\al+\int|\da^4f(\al\pm
ir)|^2d\al,
$$
and $z^0(\al)$ a curve satisfying the arc-chord condition
and $z^0(\al)\in X_{r_0}$ for some $r_0>0$, we prove the existence of a time
$T>0$ and $0<r<r_0$ so that there is a unique solution to the Muskat problem
in $C([-T,T];X_r)$. This allows us to find solutions that do not satisfy the R-T but
shrink the strip of analyticity. We extend equation \eqref{ec1d} as follows:
$$z_t(\al+i\zeta,t)=G(z(\al+i\zeta,t)),$$ with
$$G(z)(\al+i\zeta,t)=\frac{\rho^2-\rho^1}{2\pi}\int\frac{z_1(\al+i\zeta)-z_1(\al+i\zeta-\beta)}{|z(\al+i\zeta)-
z(\al+i\zeta-\beta)|^2}(\da z(\al+i\zeta)-\da
z(\al+i\zeta-\beta))d\beta.$$ For $0< r'<r$ and the open set $O$ in $S_{r}$ given by
\begin{equation}\label{openO}
O=\{z\in X_{r}: \|z\|_{r}<R,\quad
\|F(z)\|_{L^\infty(S_r)}<R^2\},
\end{equation}
the function $G$ for $G:O\rightarrow X_{r'}$ is a continuous
mapping and there is a constant $C_R$ (depending on $R$
only) such that
\begin{equation}\label{cota}
\|G(z)\|_{r'}\leq \frac{C_R}{r-r'}\|z\|_{r},
\end{equation}
\begin{equation}\label{casiL}
 \|G(z^2)-G(z^1)\|_{r'}\leq \frac{C_R}{r-r'}\|z^2-z^1\|_{r},
\end{equation}
and
\begin{equation}\label{paraarc-chord}
\sup_{\al+i\zeta\in S_r,\beta\in\T}
|G(z)(\al+i\zeta)-G(z)(\al+i\zeta-\beta)|\leq C_R|\beta|,
\end{equation}
for $z,z^j\in O$. For
initial data $z^0\in X_{r_0}$ satisfying arc-chord, we can
find a $0<r_0'<r_0$ and a constant $R_0$ such that $\|z^0\|_{r_0'}<
R_0$ and
\begin{equation}\label{arc-chord-}
\Big|\frac{(z^0_1(\al+i\zeta)-z^0_1(\al+i\zeta-\beta))^2+(z^0_2(\al+i\zeta)-
z^0_2(\al+i\zeta-\beta))^2}{\beta^2}\Big|>\frac{1}{R_0^2},
\end{equation}
for $\al+i\zeta\in S_{r_0'}$. We take $0<r<r_0'$ and $R_0<R$ to
define the open set $O$ as in \eqref{openO}. Therefore we can use
the classical method of successive approximations:
$$
z^{n+1}(t)=z^0+\int_0^t G(z^n(s))ds,
$$
for $G:O\rightarrow X_{r'}$ and $0< r'<r$. We assume by induction
that $$\|z^k\|_r(t)< R, \qquad\mbox{ and }\qquad
\|F(z^k)\|_{L^\infty(S_r)}(t)< R$$ for $k\leq n$ and $0<t<T$ with
$T=\min(T_A,T_{CK})$ and $T_{CK}$ the time obtaining in the proofs
in \cite{Nirenberg} and \cite{Nishida}. We get $\|z^{n+1}\|_r(t)< R$ that follows using (\ref{cota},\ref{casiL}).
The time $T_A$ is to yield $\|F(z^{n+1})\|_{L^\infty(S_r)}(t)< R$.
Then, using the induction hypothesis and \eqref{paraarc-chord} we can control the quantity taking $0<T_A<(R_0^{-2}-R^{-2})(C_R^2+2R_0C_R)^{-1}$.

5. Fifth, all the results above allow us to prove that there is a non empty
set of initial data in $H^4$ satisfying the arc-chord and R-T
conditions such that the solution of the Muskat problem reaches the
unstable regime: the R-T becomes strictly negative on a non-empty
interval. We pick initial data as in 3. We apply the local-existence
result in 4 to get an analytic solution $z(\al,t)$ on $[-T,T]$. Then
we consider a time $0<\delta<T$ and a curve
$\omega_\delta^\ep(\al,t)$, solving the  Muskat problem with initial
datum $z(\al,-\delta)+\eta_\delta^\epsilon(\al)$. The function
$\eta_\delta^\epsilon$ has a small  $H^4$ norm, i.e.
$$
\|\omega_\delta^\ep(\cdot,-\delta)-z(\cdot,-\delta)\|_{H^4}=
\|\eta_\delta^\epsilon\|_{H^4}\leq \ep.
$$
The time $\delta$ is small enough so that $\omega_\delta^\ep(\al,-\delta)$ satisfies R-T: $(\rho^2-\rho^1)\da(\omega_\delta^\ep)_1(\al,-\delta)>0$. Then we apply the local existence result in 1. that $\omega_\delta^\ep(\al,t)$ becomes analytic  for some time $-\delta<t$. With 2. we assure the existence and analyticity of the solution even if $\da(\omega_\delta^\ep)_1(\al,t)\leq 0$ for some time $t$. Then, we show that both solutions are close in the $H^4$ topology as time evolves. We can apply to $\omega_\delta^\ep$ the local-existence result in 4. if it is needed. Then, with $\delta$ and $\ep$ small enough we find the desired result.

\section{Turning water waves}

In this section we prove for the water wave problem ($\rho^1=0$ and
(\ref{density},\ref{mass},\ref{vorticity},\ref{Euler})) that with
initial data given by a graph $(\al,f_0(\al))$, the interface
reaches a regime in finite time where it only can be parameterized
as $z(\al,t)=(z_1(\al,t),z_2(\al,t)),$ for $\al\in\R$, with $\da
z_1(\al,t)<0$ for  $\al\in I$, a non-empty interval. Therefore there
exists a time $t^*$ where the solution of the free boundary problem
reparameterized by $(\al,f(\al,t))$ satisfies
$\|f_\al\|_{L^\infty}(t^*)=\infty$.
\begin{thm}
There exists a non-empty open set of initial data $(\al,f_0(\al))$, with $f_0\in
H^5$, such that in finite time $t^*$ the solution of the water waves
problem ($\rho^1=0$ and
(\ref{density},\ref{mass},\ref{vorticity},\ref{Euler})) given by
$(\al,f(\al,t))$ satisfies $\|f_\al\|_{L^\infty}(t^*)=\infty$. The
solution can be continued for $t>t^*$ as $z(\al,t)$ with $\da
z_1(\al,t)<0$ for $\al\in I$, a non-empty interval.
\end{thm}

In order to prove this theorem we consider a curve $z^*(\al)\in
H^5$ with the same properties as in point 3. of previous section.
Then, we pick $z(\al,t^*)=z^*(\al)$ and $\omega(\al,t^*)=\da
z^*_1(\al)$ as a datum for the initial value problem. It is easy to
find the same properties for the velocity, since  the tangential
direction does not affect the evolution. Picking the appropriate
$c(\al,t)$ and applying the local-existence result in \cite{ADP2},
there exists a solution of the water waves problem with $z(\al,t)\in
C([t^*-\delta,t^*+\delta];H^5)$, $\omega(\al,t)\in
C([t^*-\delta,t^*+\delta];H^4)$ and $\delta>0$ small enough. Then,
the initial datum
$(z_0(\al),\omega_0(\al))=(\al,f_0(\al),\omega_0(\al))$ is given by
$(z(\al,t^*-\delta),\omega(\al,t^*-\delta))$.

\section{Muskat breakdown}
In this section we show that there exists a smooth initial data in
the stable regime for the Muskat problem  such that the solution
turns to the unstable regime and later it breaks down. The outline of the
proof is to construct a curve in the unstable regime which is
analytic except in a single point. We show that as we evolve backwards in time
the curve becomes analytic and is as close as we desired (in the
$H^k$ topology with $k$ large enough) to the curve from part 3. of
section 2.

  Here we will work in the periodic setting and will consider the
equation
\begin{equation}\label{commuskat}
\pa_t
z(\zeta,t)=\int_{w\in\Gamma_+(t)}\frac{\sin(z_1(\zeta,t)-z_1(w,t))}
{\cosh(z_2(\zeta,t)-z_2(w,t))-\cos(z_1(\zeta,t)-z_1(w,t))}
(\pa_\zeta z(\zeta,t)-\pa_\zeta z(w,t))dw,
\end{equation}
where $\zeta\in \Omega(t)$, $$\Omega(t)=\{\zeta\in \C/2k\pi\,:\ |\Im \zeta
|< h(\Re z,t)\},$$ $h(x,t)$ is a positive periodic function with
period $2\pi$ and smooth for fixed time $t$ and
$$\Gamma_{\pm}(t)=\{\zeta\in \C/2k\pi\,:\ \zeta=x+ih(x,t)\}.$$ This
equation is equivalent to (\ref{density},\ref{mass},\ref{Darcy}) for
holomorphic functions.

In order to prove the result we will need the following theorem:
\begin{thm}\label{local}
Let $h(x,t)$ be a positive, smooth and periodic function with period
$2\pi$ for fixed time $t\in[t_0-\delta,t_0]$. Let $z(x,t_0)$ be a
curve satisfying the following properties:
\begin{itemize}
\item  $z_1(x,t_0)-x$ and $z_2(x,t_0)$ are periodic with
period $2\pi$.
\item $z(\zeta,t_0)$ is real for $\zeta$ real.
\item $z(\zeta,t_0)$ is analytic in $\zeta\in\Omega(t_0)$.
\item $z(\zeta,t)\in H^k(\Gamma_{\pm}(t_0))$ with $k$ a large enough integer.
\item Complex Arc-Chord condition.
$$|\cosh(z_2(\zeta,t_0)-z_2(w,t_0))-\cos(z_1(\zeta,t_0)-z_1(w,t_0))|\geq [||\Re (\zeta-w)||+|\Im(\zeta-w)|]^2,$$
for $\zeta$, $w\in\overline{\Omega}(t_0)$, where $||x||= distance(x,2k\pi).$
\item Generalized Rayleigh-Taylor condition: $RT(\zeta,t_0)>0$, where
$$RT(\zeta,t)=\Re \left(\frac{-2\pi \pa_\zeta z_1(\zeta,t)}{(\pa_\zeta z_1(\zeta,t))^2+(\pa_\zeta z_2(\zeta,t))^2}(1+i\pa_x h(\Re\zeta,t))^{-1}\right)$$
$$+\Im\left(\left\{P.V.\int_{w\Gamma_+(t)}\frac{\sin(z_1(\zeta,t)-z_1(w,t))}
{\cosh(z_2(\zeta,t)-z_2(w,t))-\cos(z_1(\zeta,t)-z_1(w,t))}dw+i\pa_t
h(\zeta,t)\right\}\right.$$
$$\left.\times(1+i\pa_x h(\Re\zeta,t))^{-1}\right)$$
\end{itemize}
Then, for small enough $\delta$, there exist a solution for the
equation (\ref{commuskat}) in the time interval
$t\in[t_0-\delta,t_0]$ satisfying
\begin{itemize}
\item  $z_1(x,t)-x$ and $z_2(x,t)$ are periodic with
period $2\pi$.
\item $z(\zeta,t)$ is real for $\zeta$ real.
\item $z(\zeta,t)$ is analytic in $\zeta\in\Omega(t_0)$.
\item $z(\zeta,t)\in H^k(\Gamma_{\pm}(t))$ with $k$ a large enough integer.
\end{itemize}
\end{thm}

Now, let $\underline{z}(x,t)$ be the solution of  the Muskat
problem with $\underline{z}(x,0)=\underline{z}^0(x)$, where
$\underline{z}^0(x)$ is the particular initial data from part 3. of
the section 2. We shall define  this solution as the unperturbed
solution. Let us denote the  Rayleigh-Taylor function
$$\sigma_1^0(x,t)\equiv\frac{-2\pi \pa_x \underline{z}_1(x,t)}{(\pa_x
\underline{z}_1(x,t))^2+(\pa_x \underline{z}_2(x,t))^2}.$$ Notice
the minus sign in the right-hand side of the previous expression.
One can check the following properties of this Rayleigh-Taylor
function:
\begin{enumerate}
\item $\sigma^0_1(\cdot,t)$ is analytic on
$\{x+iy\,:\,x\in\T,\,|y|\leq c_b\}$ with $|\sigma_1^0(x+iy,t)|\leq
C$, for all $x+iy$ as above and for all $t\leq[0,\tau]$. \item
$\sigma_1^0(0,0)$ is real for $x\in\T$, $t\in[0,\tau]$. \item
$\sigma_1^0$ has a priori bounded $C^{k_0}$ norm as a function of
$(x,t)\in\T\times[0,\tau]$ ($k_0$ large enough). \item
$\sigma_1^0(0,0)=0$. \item $\partial_x \sigma_1^0(0,0)=0$. \item
$\partial_x^2 \sigma_1^0(0,0)=-c_2<0$. \item $\partial_t
\sigma_1^0(0,0)=c_1>0$.
\end{enumerate}

 In this setting we define the following
weight functions
\begin{eqnarray}
h(x,t)&=&A^{-1}(\tau^2-t^2)+(A^{-1}-(\tau-t))\sin^2\left(\frac{x}{2}\right)\quad\text{for
$t\in[\tau^2,\tau].$}\label{dh}\\
\hbar(x,t)&=&\frac{1}{4}\left(A^{-1}\tau^2+A^{-1}\sin\left(\frac{x}{2}\right)\right)+A^{-2}\tau
t+At\sin\left(\frac{x}{2}\right)\quad\text{
$t\in[0,\tau^2].$}\label{dhbar},\end{eqnarray} with $x\in\T$. First
we choose the parameters $A$ large enough and then $\tau$ small
enough, then one can show that
\begin{equation}\label{hi}
\sigma_1^0(x,t)+\pa_th(x,t)-A^{\frac{1}{2}}h(x,t)\geq
c\tau^2\quad\text{for $x\in\T$, $t\in[\tau^2,\tau]$}\end{equation}
and
\begin{equation}\label{hbari}
\sigma_1^0(x,t)+\pa_t\hbar(x,t)-A^\frac{1}{2}\hbar(x,t)\geq
\frac{1}{2}A^{-2}\tau\quad \text{for $x\in\T$,
$t\in[0,\tau^2]$}.\end{equation} The inequalities (\ref{hi}) and
(\ref{hbari}) are one of the main ingredients of the proof of the
following results
\begin{thm}\label{mean}
Let $z(x,t)$ be a solution of the Muskat equation in the interval
$t\in[0,\tau]$. Let $h(x,t)$ and $\hbar(x,t)$ as in the expressions
(\ref{dh}) and (\ref{dhbar}) and $k$ a large enough integer. Assume
that $z(x,t)$ satisfies
\begin{itemize}
\item  $z_1(x,t)-x$ and $z_2(x,t)$ are periodic with period
$2\pi$. \item $z(\zeta,t)$ is real for $\zeta$ real. \item
$z(\zeta,t)$ is analytic in $\zeta\in\Omega(t)$. \item
$z(\zeta,t)\in H^k(\Gamma_{\pm}(t))$ with $k$ a large enough
integer. \item Complex Arc-Chord condition.
$$|\cosh(z_2(\zeta,t)-z_2(w,t))-\cos(z_1(\zeta,t)-z_1(w,t))|\geq [||\Re (\zeta-w)||+|\Im(\zeta-w)|]^2,$$
for $\zeta$, $w\in\overline{\Omega}(t)$.
\end{itemize}
Here in the definition of $\Omega(t)$ and $\Gamma_{\pm}(t)$ we use
$h(x,t)$ if $t\in[\tau^2,\tau]$ and $\hbar(x,t)$ if
$t\in[0,\tau^2].$ Then
$$\frac{1}{2}\frac{d}{dt}\left(\int_{w\in\Gamma_+(t)}
\left|\pa_\zeta^k z(\zeta,t)-\pa_\zeta^k\underline{z}(\zeta,t)\right|^2d\Re
\zeta\right)\geq -C(A)\lambda^2,$$ if $t\in [\tau^2,\tau]$
$$\int_{w\in\Gamma_+(t)}
\left|\pa_\zeta^k z(\zeta,t)-\pa_\zeta^k\underline{z}(\zeta,t)\right|^2d\Re\zeta\leq
\lambda^2$$ and $\lambda\leq \tau^{50}$.

 In addition
$$\frac{1}{2}\frac{d}{dt}\left(\int_{w\in\Gamma_+(t)}
\left|\pa_\zeta^k z(\zeta,t)-\pa_\zeta^k\underline{z}(\zeta,t)\right|^2d\Re
\zeta\right)\geq -C(A)\tau^{-1}\lambda^2,$$ if $t\in [0,\tau^2]$
$$\int_{w\in\Gamma_+(t)}
\left|\pa_\zeta^k z(\zeta,t)-\pa_\zeta^k \underline{z}(\zeta,t)\right|^2d\Re\zeta\leq
\lambda^2$$ and $\lambda\leq \tau^{50}$.
\end{thm}
This theorem implies that for all $\gamma>0$ there is $\ep>0$ such
that $$\int_{w\in\Gamma_+(t)}
\left|\pa_\zeta^k z(\zeta,t)-\pa_\zeta^k\underline{z}(\zeta,t)\right|^2d\Re\zeta\leq
\gamma$$ for $t\in[0,\tau]$ if $$\int_{w\in\Gamma_+(t)}
\left|\pa_\zeta^k z(\zeta,\tau)-\pa_\zeta^k\underline{z}(\zeta,\tau)\right|^2d\Re\zeta\leq
\ep$$ and $z(x,t)$ satisfies the requirements of the theorem.
\begin{lemma}\label{general}
Let $z(x,t)$ be  a solution of the Muskat problem satisfying the
requirements of theorem (\ref{mean}) and close enough to the
unperturbed solution in $t\in[0,\tau]$. Let $h(x,t)$ and
$\hbar(x,t)$ be as in (\ref{dh}) and (\ref{dhbar}) with a suitable
choice of $A$ and $\tau$. Then $z(x,t)$ satisfies the generalized
Rayleigh-Taylor condition in $t\in[0,\tau]$. In particular the
unperturbed solution satisfies the generalized Rayleigh-Taylor
condition in $t\in[0,\tau]$.
\end{lemma}
Theorems (\ref{local}), (\ref{mean}) and lemma (\ref{general}) allow
us to achieve the desired result. Indeed we can choose a curve
$z(x,\tau)$ such that
$$\int_{\zeta\in \Gamma_\pm}\left|\pa_\zeta^k z(\zeta,\tau)-\pa_\zeta^k \underline{z}(\zeta,\tau)\right|^2d\Re\zeta\leq
\ep,$$ with $0<\ep<\ep_0$ ($\ep_0$ small enough), satisfying the
generalized Rayleigh-Taylor condition by lemma (\ref{general}) and
satisfying the rest of the hypothesis of theorem (\ref{local}).
Since $h(0,\tau)=0$, $z(x,t)$ is allow to be no analytic at $x=0$
(maybe $z(x,\tau)\in H^k(\T)$ but $z(x,\tau)\not\in H^{k+1}(\T)$) .
By theorem (\ref{local}) there is a solution $z(x,t)$, analytic in
$\Omega(t)$, for some interval $t\in[\tau-\delta,\tau]$ with small
enough $\delta$ and for all $\ep$. By theorem (\ref{mean}), we can
choose $\ep$ small enough in such a way that, by lemma
(\ref{general}), $z(x,\tau-\delta)$ satisfies the generalized
Rayleigh-Taylor condition. Then we can go further the time
$\tau-\delta$. Iterating this argument, we find we can extend
$z(x,t)$ to be a solution of the Muskat problem, analytic in
$\Omega(t)$ for all $t\in[0,\tau]$ and  as close as we want to the
unperturbed solution.

\subsection*{{\bf Acknowledgements}}

\smallskip

 AC, DC and FG were partially supported by the grant {\sc MTM2008-03754} of the MCINN (Spain) and
the grant StG-203138CDSIF  of the ERC. CF was partially supported by
NSF grant DMS-0901040 and ONR grant ONR00014-08-1-0678. FG was partially supported by NSF grant DMS-0901810. MLF was partially supported by the grants {\sc
MTM2008-03541} and {\sc MTM2010-19510} of the MCINN (Spain).

\begin{tabular}{ll}
\textbf{Angel Castro} &  \\
{\small Instituto de Ciencias Matem\'aticas} & \\
{\small Consejo Superior de Investigaciones Cient\'ificas} &\\
{\small Serrano 123, 28006 Madrid, Spain} & \\
{\small Email: angel\underline{  }castro@icmat.es} & \\
   & \\
\textbf{Diego C\'ordoba} &  \textbf{Charles Fefferman}\\
{\small Instituto de Ciencias Matem\'aticas} & {\small Department of Mathematics}\\
{\small Consejo Superior de Investigaciones Cient\'ificas} & {\small Princeton University}\\
{\small Serrano 123, 28006 Madrid, Spain} & {\small 1102 Fine Hall, Washington Rd, }\\
{\small Email: dcg@icmat.es} & {\small Princeton, NJ 08544, USA}\\
 & {\small Email: cf@math.princeton.edu}\\
 & \\
\textbf{Francisco Gancedo} &  \textbf{Mar\'ia L\'opez-Fern\'andez}\\
{\small Department of Mathematics} & {\small Institut f\"ur Mathematik}\\
{\small University of Chicago} & {\small Universit\"at Z\"urich}\\
{\small 5734 University Avenue,} & {\small Winterthurerstr. 190}\\
{\small Chicago, IL 60637, USA}  & {\small CH-8057 Z\"urich, Switzerland}\\
{\small Email: fgancedo@math.uchicago.edu} & {\small Email:
maria.lopez@math.uzh.ch}
\end{tabular}

\end{document}